\documentclass[12pt]{article}
\usepackage{amsmath}
\usepackage{amsfonts}
\def\be{\begin{equation}}
\def\ee{\end{equation}}
\def\N{\mathbb N}
\def\R{\mathbb R}
\def\Z{\mathbb Z}
\def\qed{\hfill$\diamondsuit$}
\def\ed{\end{document}}
\def\II{I\negthinspace I}

\newtheorem{prop}{Proposition}
\newtheorem{thm}{Theorem}
\newtheorem{cor}{Corollary}
\newtheorem{lem}{Lemma}

\def\eqalign#1{\null\,\vcenter{\openup\jot\ialign
              {\strut\hfil$\displaystyle{##}$&$\displaystyle{{}##}$
               \hfil\crcr#1\crcr}}\,}
               
\title{Quenched Central Limit Theorems \\ For Sums of Stationary Processes}
\author{Dalibor Voln\'y and Michael Woodroofe\\
	The University of Rouen and The University of Michigan}

\date{}

\begin{document}
\maketitle

\section{Introduction}

	Let $(\Omega,{\cal A},\mu)$ denote a probability space and $T:\Omega \to \Omega$ an invertible, ergodic measure preserving transformation.  Thus, $T$ has a measurable inverse, $\mu[T^{-1}(A)] = \mu(A)$ for all $A \in {\cal A}$, and the only $A \in {\cal A}$ for which $T^{-1}(A) = A$ have measure $0$ or $1$.   If $f:\Omega \to \R$ is any measurable function, then $X_i = f\circ T^i,\ i \in \Z$, defines a strictly stationary process.  Conversely, any strictly stationary ergodic  process can be represented in this form on a sequence space.  Suppose that the process $X_i = f\circ T^i$ is adapted to a filtration $({\cal F}_i)_{i\in\Z}$ for which ${\cal F}_{i+1} = T^{-1}{\cal F}_i$ for all $i \in \Z$.  This is always the case with ${\cal F}_i = \sigma\{X_j: j \le i\}$, but other choices may be more convenient in other examples.  Finally, let $\mu(\cdot;\cdot)$  denote a regular conditional probability on ${\cal B} = \sigma\{X_1,X_2,\cdots\}$ given ${\cal F}_0$, \cite{L63}, p. $358-64$, and write $\mu_{\omega} = \mu(\omega;\cdot)$;  thus, $\mu(\omega;\cdot)$ is a probability measure for each $\omega \in \Omega$ and $\mu(\cdot;B)$ is (a version of) $\mu\{B|{\cal F}_0\}$ for each $B \in {\cal B}$.
		
	Next, let $U^ig = g\circ T^i$ for $i \in \Z$ and $g \in L^1(\mu)$, so that $X_i = U^if$, and let $S_n(f) = X_1+\cdots+X_{n} = Uf+\cdots+U^nf$ for $n \ge 1$.  If $f \in L_0^{2}(\mu) = \{g \in L^2(\mu): \int_{\Omega} gd\mu = 0\}$, so that $(X_i)_{i\in\Z}$ is a centered process with finite variance, then it is natural to ask whether $S_n(f)/\sqrt{n}$ is asymptotically normal as $n \to \infty$.  The question has several forms.  Let $G_n$ denote the distribution function of $S_n(f)/\sqrt{n}$, $G_n(z) = \mu\{S_n(f)/\sqrt{n} \le z\}$ for $z \in \R$, and write $G_n(\omega;z) = \mu_{\omega}\{S_n(f)/\sqrt{n} \le z\}$ for the conditional distribution function given ${\cal F}_0$.  If $G_n$ converges to a  distribution $G$, say, denoted $G_n \Rightarrow G$, then the convergence is said to be {\it annealed}.  If there is single (non-random) distribution $G$ for which $G_n(\omega;\cdot) \Rightarrow G$ for $a.e.\ \omega$, then the convergence is said to be {\it quenched} with respect to ${\cal F}_0$ or simply ${\cal F}_0$-quenched (Derriennic and Lin \cite{DL01}).  There is another possibility in which almost everywhere convergence is replaced by convergence in probability.  We will call such convergence {\it weakly ${\cal F}_0$-quenched}.   Necessary and sufficient conditions for weakly quenched convergence may be found in \cite{DM02} and \cite{WW04}.  Conditions for fully quenched convergence are more delicate, for example \cite{C09}, \cite{CP09} and \cite{ZW08}.  

	Interest in quenched limits comes in part from additive functionals of a Markov Chain that does not start in the stationary distribution.  For example, in Markov Chain Monte Carlo, the stationary distribution $\pi$ contains an unknown normalizing constant, and the purpose is to generate random variables having approximately the distribution $\pi$.  Thus, let $\cdots W_{-1},W_0,W_1,\cdots$ denote a stationary ergodic Markov Chain with stationary marginal distribution $\pi$ and consider an additive functional $S_n = g(W_1)+\cdots+g(W_n)$ where $g \in L^2(\pi)$.  This is a stationary process and, so, can be represented in the form $X_k = f\circ T^k$ as above.  If ${\cal F}_k = \sigma\{W_i: i \le k\}$, then the question of quenched convergence is equivalent to asking whether the conditional distribution of $S_n/\sqrt{n}$ given $W_0=w$ converges for $a.e.\ w\ (\pi)$.

	Denote the norm in $L^p(\mu)$ by $\Vert\cdot\Vert_p$; write $E$ for both expectation and conditional expectation with respect to $\mu$; and let ${\cal M}$ be the set of ${\cal F}_0$ measurable  $g \in L^2(\pi)$ for which  $E(Ug|{\cal F}_{0}) = 0$.  ${\cal M}$ is called the {\it martingale difference space} because $U^ig, i \in \Z$, are martingale differences for any $g \in {\cal M}$.  An $f \in L^2(\mu)$ is said to admit an $L^p$ co-boundary if there is an $m \in {\cal M}$ and a $g \in L^p(\mu)$ for which $f = m + g-Ug$.  It is known that the existence of an $L^1$ co-boundary implies the annealed version the the Central Limit Theorem and the existence of an $L^2$ co-boundary the quenched version.  In Theorem 2 we show that the existence of an $L^1$ co-boundary does not imply the quenched version.  In Theorem 1 it is shown that the condition (\ref{eq:hnn}) introduced by Hannan \cite{H73} implies that the conditional distributions of $\{S_n-E(S_n|{\cal F}_0)\}/\sqrt{n}$ converge to a normal distribution $w.p.1$.

\section{The Martingale Case}

	Much of the recent progress on the central limit question for sums of stationary processes has relied on approximation by martingales.  So, the case in which $S_n(f)$ is a martingale is considered first.   Some preparation is necessary.  Recall that if $g \in L^1(\mu)$, then $U^ng/n \to 0\ w.p.1$ by the Pointwise Ergodic Theorem (since $U^ng = S_n(g) -S_{n-1}(g)$).  So, if $g \in L^2(\mu)$, then $U^ng/\sqrt{n} \to 0\ w.p.1$.  Next, recall Mc Leish's conditions for the Martingale Central Limit Theorem \cite{M74}, as refined by Lachout \cite{L85}:  Let $X_{n,j}, 1 \le j \le n$, be martingale differences  defined on a probability space $(\Omega_n,{\cal A}_n,P_n)$ for each $n$.   If

\smallskip

	(i) $\sum_{j=1}^n X_{n,j}^2 \to^{P_n} \sigma^2$,
	
\smallskip

	(ii) $\forall\ \epsilon > 0,\ P_n[\max_{j\le n} |X_{n,j}| \ge \epsilon] \to 0$,
	
\smallskip

	(iii) $\sup_{n\ge 1} E_n[\max_{j\le n} X_{n,j}^{2}] < \infty$,
 
\smallskip\noindent
then $X_{n,1}+\cdots+X_{n,n} \Rightarrow {\rm Normal}[0,\sigma^2]$ (that is, the distribution function of the sum converges weakly to the normal distribution with mean $0$ and variance $\sigma^2$.)  Lachout \cite{L85} showed that (ii) and (iii) could be replaced by

\smallskip

	(iv)  $E_n[\max_{j \le n} |X_{n,j}|] \to 0$.

\smallskip\noindent
	
\begin{lem}\label{lem:iv}
	If $X_{n,j} = f\circ T^j/\sqrt{n}$, where $f \in L^2(\mu)$, then (iv) holds.
\end{lem}	

	{\it Proof}. Since $f \in L^2(P)$, it is clear that $\lim_{n\to\infty} U^nf/\sqrt{n} = 0\ w.p.1\ (\mu)$ and, therefore, that $\lim_{n\to\infty} \max_{j\le n} |U^jf|/\sqrt{n} = 0\ w.p.1.\ (\mu_{\omega})$ for $a.e.\ \omega$.  Let 
$$
	X^* = \sup_{n\ge 1} \max_{1\le j\le n} |X_{n,j}|\quad {\rm and}\quad  f^* = \sup_{n\ge 1} \sqrt{{f^2\circ T+\cdots+f^2\circ T^n\over n}}.
$$
Then $X^* \le f^*$ clearly, and $\mu\{f^* > \lambda\} \le \Vert f\Vert_2^{2}/\lambda^2$ for $\lambda > 0$, by the Maximal Ergodic Theorem, \cite{B95}, p. $318$.  So, $f^* \in L^1(\mu)$ and, therefore, $X^* \in L^1(\mu_{\omega})$ for $a.e.\ \omega$. Condition  (iv) then follows from the Dominated Convergence Theorem, applied conditionally. \qed

\begin{prop}\label{prop:mrtngl}
	If $f \in {\cal M}$, then $G_n(\omega;\cdot)$ converges to the normal distribution with mean $0$ and variance $\Vert f\Vert_2^{2}$ for $a.e.\ \omega$.
\end{prop}

	{\it Proof}.  Equivalent results are established in \cite{KV86}, \cite{W92}, and \cite{WW04}, but are not isolated there.  A proof is included here for completeness.  Our proof, at least, differs from the arguments of \cite{KV86}, \cite{W92}, and \cite{WW04}.   Let $X_{n,j} = U^jf/\sqrt{n},\ 1 \le j \le n$.  Then $X_{n,j}$ are martingale differences with respect to $\mu_{\omega}$ for $a.e.\ \omega$, by the smoothing property of conditional expectation, \cite{B95}, p. $448$.   So, it suffices to show that (i) and (iv) hold with $P_n = \mu_{\omega}$ for $a.e.\ \omega$.  That (iv) is satisfied was shown in Lemma \ref{lem:iv}.  For (i), it follows directly from the Pointwise Ergodic Theorem that
$$
	\sum_{j=1}^n X_{n,j}^{2} ={1\over n}\sum_{j=1}^n U^jf^2 \to \Vert f\Vert_2^{2},
$$
$w.p.1\ (\mu)$ and, therefore, $w.p.1\ (\mu_{\omega})$ for $a.e.\ \omega\ (\mu)$.  So, (i) is satisfied for the measures $\mu_{\omega}$ for $a.e.\ \omega$.   \qed

\begin{cor}
	If $f$ admits a $L^2$-coboundary, say $f = m + g-Ug$, where $m \in {\cal M}$ and and $g \in L^2(\mu)$, then $G_n(\omega;\cdot) \Rightarrow {\rm Normal}[0,\Vert m\Vert_2^{2}]$ for $a.e.$.
\end{cor}

	{\it Proof}. Since $S_n(f) = S_n(m) + Ug-U^{n+1}g = S_n(m)+o(\sqrt{n})\ w.p.1$, the corollary follows from Slutzky's Theorem, applied conditionally. \qed

\section{Hannan's Theorem}

	Write $E^i(g) = E(g|{\cal F}_i)$ and $P_ig = E^i(g)-E^{i-1}(g)$ for $g \in L^1(\mu)$ and observe that the restrictions of each $E^i$ and $P_i$ are projections on $L^2(\mu)$.  Hannan \cite{H73} has shown that if $f \in L^2(\mu)$ is ${\cal F}_{0}$ measurable, $E(f|{\cal F}_{-\infty}) = 0$, and 
\be\label{eq:hnn}
	\sum_{i=-\infty}^0 \Vert P_if\Vert_2 < \infty, 
\ee
then $S_n(f)/\sqrt{n} \Rightarrow {\rm Normal}[0,\Vert m\Vert_2^{2}]$, where $m = \sum_{i=0}^{\infty} P_0U^if$.    Voln\'y and Woodroofe \cite{VW09} show by example that the convergence need not be quenched.  However,

\begin{thm}\label{thm:hnn}
	If $f \in L_0^{2}(\mu)$ is ${\cal F}_0$-measurable, $E(f|{\cal F}_{-\infty}) = 0$ and $\sum_{i\le 0} \Vert P_if\Vert_2 < \infty$,  then the conditional distribution of $\{S_n(f)-E^0[S_n(f)]\}/\sqrt{n}$ given ${\cal F}_0$ converges to a normal distribution $w.p.1$.
\end{thm}

	{\it Proof}.  The easily verified relations $E^i(U^jf) = U^jE_{i-j}f$ and $P_i(U^jf) = U^jP_{i-j}f$ are used in the proof.  Letting $g_n = P_0U^{-1}S_n(f)$,
$$
	\eqalign{{S_n(f)-E^0[S_n(f)]} &= \sum_{j=1}^n P_j[S_n(f)-S_{j-1}(f)]\cr 
					&= \sum_{j=1}^n P_jU^{j-1}S_{n-j+1}(f) = \sum_{j=1}^n U^jg_{n-j+1}.\cr}
$$
So, 
$$
	{S_n(f)-E^0[S_n(f)]\over\sqrt{n}} = \sum_{j=1}^n Y_{n,j},\ 
$$	
where $Y_{n,j} = U^jg_{n-j+1}/\sqrt{n}$.  The $Y_{n,j}$ are martingale differences, since $E^{j-1}[U^jg_{n-j+1}] = U^jE^{-1}g_{n-j+1} = 0$.  So, it suffices to show that conditions (i) and (iv) of Proposition \ref{prop:mrtngl} hold $w.p.1\ (\mu)$.  That (iv) is satisfied follows from Lemma \ref{lem:iv}. For (i) first observe that $|g_n| \le \sum_{k=1}^{\infty} |P_0U^kf| = h$, say, and $\Vert h\Vert_2 \le \sum_{k=1}^{\infty} \Vert P_0U^kf\Vert = \sum_{k=1}^{\infty} \Vert P_{-k}f\Vert$, which is finite by assumption.  So, $g_n$ converges to $m = \sum_{k=1}^{\infty} P_0U^kf\ w.p.1$ and in $L^2(P)$, $E[\sup_{n\ge 1} |g_n|^2] < \infty$, and 
$$
	\sum_{j=1}^n Y_{n,j}^2 = {1\over n}\sum_{j=1}^n (U^jg_{n-j+1})^2 \to \Vert m\Vert_2^{2}\ w.p.1\ (\mu)
$$
by (a simple application of) the The Pointwise Ergodic Theorem.  See Lemma \ref{lem:smb} below.  \qed

\begin{lem}\label{lem:smb}
	If $g_n$ are measurable, $g_n \to g\ a.e.$, and $\sup_{n\ge 1} |g_n|$ is integrable, then
$$
	{1\over n}\sum_{k=1}^n U^kg_{n-k+1} \to \int_{\Omega} gdP\ w.p.1.
$$
\end{lem}

	{\it Proof}.  There is no loss of generality in supposing that $g = 0$.  Let $h_m = \sup_{n\ge m} |g_n|$.  Then $U^nh_1 = o(n)\ w.p.1$, since $h_1$ is integrable, and
$$
	\left| {1\over n}\sum_{k=1}^n U^kg_{n-k+1}\right| \le {1\over n}\sum_{k=1}^{n-m} U^kh_m + {1\over n}\sum_{k=1}^m U^{n-k+1}h_1 \to \int_{\Omega} h_mdP
$$
for any $m \ge 1$.  The lemma follows since the right side may be made arbitrarily small by taking $m$ sufficiently large. \qed

	A result equivalent to Theorem \ref{thm:hnn} appears in \cite{CP09}; the two results were obtained independently.
	
\section{Co-boundaries}

	In this section it is shown that the existence of an $L^1$ co-boundary does not imply the quenched version of the central limit theorem.  The same construction shows that a condition suggested by Heyde \cite{H75} does not imply quenched convergence.

\begin{thm}\label{thm:cbndry}
	There is a dynamical system $(\Omega,{\cal A},\mu,T)$, an $f \in L_0^{2}(\mu)$, and a $g \in L^1(\mu)$ for which:

\smallskip
	(a) $f = g-Ug$;
	
\smallskip
	(b) $\limsup_{n\to\infty} E^0[S_n(f)]/\sqrt{n} \ge 1\ w.p.1$;
		
\end{thm}

	{\it Proof}.  Let $\ell_1,\ell_2, \cdots$ and $M_1,M_2,\cdots$ be strictly increasing sequences of positive integers for which 
\be\label{eq:cndtn1}
	\sum_{k=1}^{\infty} \left[{1\over\sqrt{k\ell_k}} +  \ell_k^{2}\sqrt{M_{k-1}\over M_{k}}\right]  < \infty.
\ee
If (\ref{eq:cndtn1}) holds for $M_k$ then it also holds with $M_k$ replaced by $M_k\vee 2\ell_k$, and we will suppose that $M_k \ge 2\ell_k$ for all $k$.  Let $e_k$ be random variables for which $U^ie_k = e_k\circ T^i, i \in \Z,\ k \in \N$ are independent, 
$$
	e_k = \pm {\sqrt{M_k}\over \ell_k}\ {\rm with\ probability}\  {1\over 2kM_k}
$$ 
each, and $e_k = 0$ otherwise; and let ${\cal F}_j = \sigma\{U^ie_k: i \le j,\ k \ge 1\}$.   Observe that $\Vert e_k\Vert_1 = 1/k\ell_k\sqrt{M_k}$ and $\Vert e_k\Vert_2 ^{2} = 1/k\ell_k^{2}$ for $k \ge 1$.  If $\ell_k$ is any sequence for which $1/\sqrt{k\ell_k}$ summable, then (\ref{eq:cndtn1}) holds with $M_0 = 1$ and $M_k = k\ell_k^{5}M_{k-1}$ for $k \ge 1$.

	To construct $f$ and $g$ and verify (a), first let $N_k = M_1+\cdots+M_k$ and observe that $N_{k-1} \le kM_{k-1} = o(M_k)$ as $k \to \infty$ by (\ref{eq:cndtn1}), so that $N_k \sim M_k$.  Next, let $h_k =  e_k+\cdots+U^{\ell_k-1}e_k$, and 
\be\label{eq:eff}
 	f_k = h_k - U^{\ell_k}h_k = \sum_{i=0}^{\ell_k-1} U^ie_k - \sum_{i=\ell_k}^{2\ell_k-1} U^{i}e_k.
\ee
Then $\Vert f_k\Vert_2^{2} = 2\ell_k\Vert e_k\Vert^2 = {2/ k\ell_k}$.  So, $\Vert f_k\Vert_2$ is summable by (\ref{eq:cndtn1}) and, therefore, $f = \sum_{k=1}^{\infty} U^{-N_k}f_k \in L^2(\mu)$.  Next, let $g_k = h_k+\cdots+U^{\ell_k-1}h_k$.  Then $g_k - Ug_k = h_k - U^{\ell_k}h_k = f_k$, and  $\Vert g_k\Vert_1 \le \ell_k\Vert h_k\Vert_1 \le \ell_k^{2}\Vert e_k\Vert_1 = {\ell_k/ k\sqrt{M_k}}$, which is also summable by (\ref{eq:cndtn1}).  So, $g = \sum_{k=1}^{\infty} U^{-N_k}g_k \in L^1(\mu)$ and $f = g - Ug$, establishing (a).  For later reference, observe that $S_n(f) = Ug-U^{n+1}g$ and 
\be\label{eq:geek}
	g_k = \sum_{i=0}^{\ell_k-1}\sum_{i=0}^{\ell_k-1} U^{i+j}e_k = \sum_{r=0}^{2\ell_k-2} c_{k,r}U^re_k,
\ee
where $c_{k,r} = (r+1)$ for $0 \le r \le \ell_k-1$ and $c_{k,r} = (2\ell_k-1-r)$ for $\ell_k \le r \le 2\ell_k-2$.

		To verify (b) first observe that $E^0(U^{n}g) =  I(n) + \II(n)$,  where 
$$
	I(n) = \begin{cases}U^{n-N_k}g_k&{\rm if}\ N_{k-1} < n \le N_k -2\ell_k + 2 \\ \sum_{i=0}^{N_k-n} c_{k,i}U^{i+n-N_k}e_k&{\rm if}\ N_k -2\ell_k + 2 <  n \le N_k
			\end{cases}
$$
and
$$
	\II(n) = \sum_{j=k+1}^{\infty} U^{n-N_j}g_j
$$
for $N_{k-1} < n \le N_k$.  It will be shown that
\be\label{eq:one}
	P\left[\max_{N_{k-1} < n \le N_k} I(n) \ge (1-\epsilon)\sqrt{N_k},\ {\rm i.o.} \right] = 1
\ee
for any $\epsilon > 0$ and
\be\label{eq:two}
	P\left[\max_{N_{k-1} < n \le N_k}\II(n) \ne 0\ {\rm i.o}\right] = 0.
\ee
Assertion (b) then follows from $E^0[S_n(f)] = Ug - E^0[U^{n+1}g]$.  

	To establish (\ref{eq:one}),  let $A_{n}$ be the event that $U^{n+\ell_k-1-N_k}e_k = {\sqrt{M_k}/\ell_k}$ and $U^{n+i-N_k}e_k = 0\ {\rm for}\ 0 \le i\ne\ell_k-1\le 2\ell_k-2$ for $N_{k-1} < n \le N_k$.  If $\epsilon > 0$, then 
$$
	\bigcup_{n=N_{k-1}+1}^{N_k-2\ell_k+2} A_{n} \subseteq \left\{\max_{N_{k-1}<n\le N_k} I_k(n) \ge (1-\epsilon)\sqrt{N_k}\right\},
$$
for all large $k$, since $N_k \sim M_k$.  Then from (\ref{eq:cndtn1}) and the Bonferoni Inequalities, \cite{F68}, p. $100$
$$
	\eqalign{P\left(\bigcup_{n=N_{k-1}+1}^{N_k-2\ell_k+2} A_{n}\right) &\ge \sum_{n=N_{k-1}+2}^{N_k-2\ell_k+2} P(A_{n}) - \sum_{n=N_{k-1}+1}^{N_k-2\ell_k+2} \sum_{m=N_{k-1}+1}^{n-1} P(A_{m}\cap A_{n})\cr
			&\ge (N_k-N_{k-1}-2\ell_k+2)\left({1\over 2kM_k}\right)\left[1- {1\over kM_k}\right]^{2\ell_k-2}\cr 
			&- {1\over 2}(N_k-N_{k-1}-2\ell_k)^2\left({1\over 2kM_k}\right)^2\left[1- {1\over kM_k}\right]^{2\ell_k-2}\cr
			&\ge {1\over 2k} + o({1\over k}).\cr}
$$
That (\ref{eq:one}) holds, then follows from the Borel Cantelli Lemmas, since the events $\cup_{n=N_{k-1}+1}^{N_k-2\ell_k+2} A_{n}$ are independent.

	For (\ref{eq:two}) simply observe that
$$
	P\left[\max_{N_{k-1}<n\le N_k} |U^{n-N_j}g_j| \ge 0\right] \le \sum_{i=N_{k-1}-N_j+1}^{N_k-N_j+2\ell_j} P\left[U^ie_j \ne 0\right]  = {N_k-N_{k-1}+2\ell_j\over jM_j}
$$
for $j > k$ and
$$
	\sum_{k=1}^{\infty}\sum_{j=k+1}^{\infty} \left({N_k-N_{k-1}+2\ell_j\over jM_j}\right) \le \sum_{j=2}^{\infty}  {N_{j-1}+2j\ell_j\over jM_j},
$$
which is finite (\ref{eq:cndtn1}).  So, (\ref{eq:two}) also follows from the Borel Cantelli Lemmas. \qed

	Recall that ${\cal M}$ denotes the martingale difference space.  Recall too the Convergence of Types Theorem, \cite{L63}, p. $203$:  {\it Let $a_n > 0$, $b_n \in \R$, and let $Y_n$ be a random variables.  If $Y_n \Rightarrow Y$, and $a_nY_n+b_n \Rightarrow Z$, where $Y$ has a non-degenerate distribution, then $a = \lim_{n\to\infty} a_n$ and $b = \lim_{n\to\infty} b_n$ exist, and $Z =^{Dist} aY+b$}.
	
\begin{cor}\label{cor:em}
	If $m \in {\cal M}$, then $S_n(m+f)/\sqrt{n} \Rightarrow {\rm Normal}[0,\Vert m\Vert_2^{2}]$, but the convergence is not quenched.
\end{cor}

	{\it Proof}.  Write
$$
	\eqalign{{S_n(m+f)\over\sqrt{n}} &= {S_n(m)\over\sqrt{n}} + {S_n(f)-E^0(S_n(f))\over\sqrt{n}} + {E^0(S_n(f))\over\sqrt{n}}\cr 
								&= Y_{n,1}+Y_{n,2}+\nu_n,\cr}
$$
say, and observe that $\nu_n$ is ${\cal F}_0$-measurable.   By Proposition \ref{prop:mrtngl} the conditional distribution of $Y_{n,1}$ given ${\cal F}_0$ converges to Normal$[0,\Vert m\Vert^2]$.  In this case$ Y_{n,2}\{S_n(f)-E^0[S_n(f)]\}/\sqrt{n}$ and ${\cal F}_0$ are independent.  So, the conditional distribution of $Y_{n,2}$ is the same as its unconditional distribution which converges to the degenerate distribution at $0$ by part (a) of the Theorem.  It follows that the conditional distribution of $Y_n = Y_{n,1}+Y_{n,2}$ converges to the normal distribution with mean $0$ and variance $\Vert m\Vert_2^{2}$.  That the conditional distribution of $S_n(m+f)/\sqrt{n}$ does not converge to Normal$[0,\Vert m\Vert^2]$ then follows from the Convergence of Types Theorem in the case $m \ne 0$. If $m = 0$, then $G_n(\omega;z) = F_n[z-n^{-{1\over 2}}E^0(S)(\omega)]$, where $F_n$ denotes the uncondtional distribution of $Y_n$.  It then follows from (b) that $\liminf_{n\to\infty} G_n(\omega;z) = 0$ for $0 < z < 1/2$. \qed 

	Heyde \cite{H75} showed that the CLT is true under a slightly weaker condition than Hannan's,
\be\label{eq:hyd}
	m = \sum_{i\in\Z} P_0U^if\quad{\rm and}\quad \Vert m\Vert_2^{2} = \lim_{n\to\infty} {E[S_n(f)^2]\over n}.
\ee
The last corollary shows that the $f$ constructed in Theorem \ref{thm:cbndry} satisfies (\ref{eq:hyd}), so that (\ref{eq:hyd}) does not imply quenched convergence either.  That (\ref{eq:hyd}) does not imply the weak invariance principle was shown in \cite{DMV07}
	
\begin{cor}\label{cor:heyde}
	If $\ell_k = 2^k$, then the $f$ constructed in Theorem \ref{thm:cbndry} satisfies (\ref{eq:hyd}) with $m = 0$.
\end{cor}

	{\it Proof}.  Since $f$ is ${\cal F}_0$-measurable, the sum on the left side of (\ref{eq:hyd}) is $\lim_{n\to\infty} P_0S_n(f)$.  If $n \ge 1$ then, using (\ref{eq:geek}),
$$
	P_0U^{n-N_k}g_k = \begin{cases} c_{k,N_k-n}e_k&{\rm if}\ N_k-2\ell_k+2 < n \le N_k \\
							0&{\rm otherwise}
						\end{cases}
$$
Let $k = k_n$ be the unique integer for which $N_{k-1} < n \le N_k$. Then $S_n(f) = P_0U^{n-N_k}g_{N_k}$ and $\Vert P_0S_n(f)\Vert_2^{2} \le 2\ell_k^{2}\Vert e_k\Vert_2^{2} = 2/k \to 0$ as $n \to \infty$.  Thus the sum converges to $m = 0$.

	For the second part of (\ref{eq:hyd}), observe that
$$
	g_k-U^ng_k = \sum_{r=0}^{n\wedge\ell_k-1} U^rh_k - \sum_{r= n\vee\ell_k-1}^{n+\ell_k-1} U^rh_k.
$$
So,
$$
	S_n(f) = \sum_{k=1}^{\infty} U^{-N_k}\left[\sum_{r=0}^{n\wedge\ell_k-1} U^rh_k\right] - \sum_{k=1}^{\infty} U^{-N_k}\left[\sum_{r= n\vee\ell_k-1}^{n+\ell_k-1} U^rh_k\right]
$$
$$
	\Vert S_n(f)\Vert_2^{2} \le 2\sum_{k=1}^{\infty} \Big\Vert \sum_{r=0}^{n\wedge\ell_k-1} U^rh_k\Big\Vert^2 + 2\sum_{k=1}^{\infty} \Big\Vert \sum_{r= n\vee\ell_k-1}^{n+\ell_k-1} U^rh_k\Big\Vert_2^{2} 
$$
and, therefore,
$$
	\Vert S_n(f)\Vert_2^{2} \le 4\sum_{\ell_k \le n} {\ell_k^{3}\over k\ell_k^{2}} + 4\sum_{\ell_k > n} {n^2\ell_k\over k\ell_k{2}}
		= 4\sum_{\ell_k \le n} {\ell_k\over k} + 4\sum_{\ell_k > n} {n^2\over k\ell_k}.
$$
Suppose now that $\ell_k = 2^k$ and let $m_n = \log_2(n)$.  Then the last term is at most
$$
	\sum_{k\le{1\over 2}m_n} 2^k + {2\over m_n}\sum_{k\le m_n} 2^k \le 2^{{1\over 2}m_n} +{2\over m_n}2^{m_n} + n^2\sum_{k>m_n} {1\over k2^k} = o(n),
$$
as required in (\ref{eq:hyd}). \qed
	
	The proofs of Theorem \ref{thm:cbndry} and Corollary \ref{cor:heyde} have been adapted and simplified from \cite{DMV07}

\ed
	
If $H_n(\omega:\cdot) \Rightarrow G$ for $a.e.\ \omega\ (\mu)$